\documentclass[12pt]{article}
\usepackage{amsmath, amsthm, amssymb, amsxtra}
\newcommand{\Lpp}{\{L_{p_0}(U_0),\,L_{p_1}(U_1)\}}
\newcommand{\Lppp}{L_{p_0}(U_0),\,L_{p_1}(U_1)}
\newcommand{\Lqq}{\{L_{q_0}(V_0),\,L_{q_1}(V_1)\}}
\newcommand{\Lqqq}{L_{q_0}(V_0),\,L_{q_1}(V_1)}
\newcommand{\lpu}{l_{p_0},\,l_{p_1}(\uuu)}

\newcommand{\lpp}{\{l_{p_0},\,l_{p_1}(\tttt)\}}
\newcommand{\lqq}{\{l_{q_0},\,l_{q_1}(\uuu)\}}

\newcommand{\lrrr}{l_{r_0},\,l_{r_1}(\uuu)}
\newcommand{\ii}{\infty}
\newcommand{\XX}{\{X_0,\,X_1\}}
\newcommand{\YY}{\{Y_0,\,Y_1\}}
\newcommand{\ff}{\varphi}
\newcommand{\rr}{\rho}
\newcommand{\pp}{\psi}
\newcommand{\uu}{u_m}
\newcommand{\uuu}{u_m^{-1}}
\newcommand{\ttt}{t_n}
\newcommand{\tttt}{t_n^{-1}}
\newcommand{\0}{r_0^{-1}=(q_0^{-1}-p_0^{-1})_+}
\newcommand{\1}{r_1^{-1}=(q_1^{-1}-p_1^{-1})_+}

\begin{document}
\begin{center}
{\Large \bf Interpolation orbits in the Lebesgue spaces}

{\sc V.I.Ovchinnikov}
\end{center}
\hyphenation{}

This paper is devoted to description of interpolation orbits with
respect to linear
operators mapping an arbitrary couple of $L_p$ spaces with weights
$\{L_{p_0}(U_0),\,L_{p_1}(U_1)\}$ into an arbitrary couple
$\Lqq$, where $1\le p_0,p_1,q_0, q_1\le\ii$.
By $L_p(U)$ we denote the space of
measurable functions $f$ on a measure space ${\frak M}$ such that $fU\in L_p$
with
the norm $\|f\|_{L_p(U)}=\|fU\|_{L_p}$.

Let $\{X_0,\,X_1\}$ and $\{Y_0,\,Y_1\}$ be two Banach couples,
$a\in X_0+X_1$. The space $Orb(a, \XX\to\YY)$ is a Banach space of
$y\in Y_0+Y_1$
such that $y=Ta$, where $T$ is a linear operator mapping
the couple $\XX$
into the couple $\YY$. The norm is defined by
$$
\|y\|_{Orb}=\inf\limits_{T}\max (\|T\|_{X_0\to Y_0},\,\|T\|_{X_1\to Y_1}),
$$
where infimum is taken over all representations of $y$ in the form $y=Ta$.
This space is called an interpolation orbit of the element $a$.

We are going to describe the spaces
$$
Orb(a,\,\Lpp\to\Lqq)
$$
for arbitrary $a$, any $1\le p_0,\,p_1,\,q_0,\,q_1\le \ii$ and
any positive weights $U_0,\,U_1,\,V_0,\,V_1$.

Fundamental results on description of these spaces in separate cases are well
known since 1965.
The key role was played by the J.Peetre K-functional.

Let $\XX$ be a Banach
couple, $x\in X_0+X_1$, $s>0$, $t>0$.
Denote by
$$
K(s,t, x;\XX)=\inf\limits_{x=x_0+x_1}s\|x_0\|_{X_0}+t\|x_1\|_{X_1}
$$
where infimum is taken over all representations of $x$ as a sum of
$x_0\in X_0$ and $x_1\in X_1$. The function $K(s,t)$ is concave and
it is uniquely defined by the function $K(1,t, x;\XX)$ which is also denoted by
$K(t, x;\XX)$.

If $1\le p_0\le q_0\le \ii$, $1\le p_1\le q_1\le \ii$,
the orbits $Orb(a,\Lpp\to\Lqq)$ were described as the generalized
Marcinkiewicz spaces, i.e.,
$$
Orb\,(a,\Lpp\to\Lqq)
$$
$$
=M_\ff(\Lqqq)=\{y:\sup_{s,t}\frac{K(s,t,y;\Lqq)}{\ff(s,t)}<\ii\},
$$
where $\ff(s,t)=K(s,t,a;\Lpp)$, for any $a\in
L_{p_0}(U_0)+L_{p_1}(U_1)$. First results were obtained by
B.S.Mitiagin and A.P.Calderon in the case $p_0=q_0=1$ and
$p_1=q_1=\ii$.
The final steps in the description were done by
G.Sparr in \cite{S}, \cite{SS} and V.I.Dmitriev in \cite{D}. In
particular G.Sparr showed that if
$$
K(s,t,y;\{L_{p_0}(V_0),\,L_{p_1}(V_1)\})\le C K(s,t,a;\Lpp),
$$
then there exists a linear operator $T:\Lpp\to\{L_{p_0}(V_0),\,L_{p_1}(V_1)\}$
such that $y=Ta$.

V.I.Dmitriev in \cite{D} had also found
a description of orbits in the case of arbitrary $1\le p_0,\,p_1\le\ii$
and $q_0=q_1=1$ as well as in the case of arbitrary $1\le p_1,\, q_0\le\ii$
and $p_0=q_1=1$. These Dmitriev's results were not estimated
properly at that time
because the form of description seemed to be not satisfactory.

The result we are going to present here goes up to the paper \cite{4}
where some
optimal interpolation theorems were found. These theorems could be considered
as a description of interpolation orbits of elements $a$ from a
one-parameter family
but for arbitrary indices $1\le p_0,p_1,q_0, q_1\le\ii$.
Developing this approach the following description
was conjectured in \cite{5}. The space $Orb(a,\,\Lpp\to\Lqq)$ is situated between
$L_{q_0}(V_0)$ and $L_{q_1}(V_1)$ exactly in "the same place" where the
Calderon--Lozanovskii space
$\ff(L_{r_0}(W_0),\,L_{r_1}(W_1))$ is situated between
$L_{r_0}(W_0)$ and $L_{r_1}(W_1)$,
where $\ff(s,t)=K(s,t,a,\{L_{p_0}(U_0),\,L_{p_1}(U_1)\})$ and
$r_0^{-1}=(q_0^{-1}-p_0^{-1})_+$, $\1$. This hypothesis
was confirmed in the paper \cite{6} for elements $a$ with quasi-power
functions
$\ff(s,t)=K(s,t,a,\{L_{p_0}(U_0),\,L_{p_1}(U_1)\})$. (Note that it may happen
that the space $L_{p_0}(U_0)+L_{p_1}(U_1)$ does not contain such elements
at all.)
Now we show that hypothesis from \cite{5}
is true for any $a\in L_{p_0}(U_0)+L_{p_1}(U_1)$. We also present a slightly
modified description of interpolation orbits which resembles Dmitriev's
description from \cite {D}.

$1^o$. {\sc The method of means for any quasi-concave functional
 parameter.}

 Let $\ff(s,t)$ be interpolation function, that is let $\rr(t)=\ff(1,t)$
 be quasi-concave and $\ff(s,t)$ be homogeneous of the degree one.
 Assume that $\ff\in\Phi_0$ which means that $\ff(1,t)\to 0$ and
$\ff(t,1)\to 0$ as $t\to 0$. Denote by $\{\ttt\}$ a balanced
sequence corresponding to the quasi-concave function $\rr(t)$ (see
\cite{bs}). The sequence may be constructed by induction
\begin{equation}\label{J_n}
\min \left(
\frac{\rr(t_{n+1})}{\rr(\ttt)},
\,\frac{t_{n+1}\rr(\ttt)}{\ttt\rr(t_{n+1})}
\right)=
q>1.
\end{equation}
(For simplicity in the sequel we suppose that the sequence
$\{\ttt\}$ is two-sided.) Such sequences invented by K.Oskolkov
were introduced to interpolation by S.Janson (see \cite{10}).

The main property of this sequence is the following
\begin{equation}\label{J}
K(s,t,\{\rr(\ttt)\},\lpp)\asymp \ff(s,t)
\end{equation}
for any $1\le p_0,p_1\le\ii$ (see \cite{8}).

\noindent{\bf Definition 1.} \textit{Let $\XX$ be any Banach
couple, $\rr(t)$ be a quasi-concave function such that
$\ff\in\Phi_0$ and $1\le p_0,p_1\le\ii$. Denote by
$\ff(X_0,\,X_1)_{p_0,p_1}$ the space of $x\in X_0+X_1$ such that
$$
x=\sum_{n\in \Bbb Z}\rr(\ttt)w_n \quad ({\rm convergence\ in}\ X_0+X_1),
$$
where $w_n\in X_0\cap X_1$ and $\{\|w_n\|_{X_0}\}\in l_{p_0}$,
$\{\ttt\|w_n\|_{X_1}\}\in l_{p_1}$.}

The norm in $\ff(X_0,\,X_1)_{p_0,p_1}$ is naturally defined. In
the case of $\ff(s,t)=s^{1-\theta}t^\theta,\quad 0<\theta<1$ these
spaces where introduced by J.-L.Lions and J.Peetre and were called
the spaces of means (see \cite{7}). In this case we take
$t_n=2^{n}$. As we shall see later that new notations we use here
have it's own reasons.

Note that the space $\ff(X_0,\,X_1)_{\ii,\ii}$ coincides with the generalized
Marcinkiewicz space $M_{\ff}(X_0,\,X_1)$ as well as with the space
$(X_0,\,X_1)_{\rr,\ii}$ (see \cite{10},\cite{8}).

Let $\XX$ be a couple of Banach lattices. Recall that the space
$\ff(X_0,\,X_1)$
 is the space of all elements from $X_0+X_1$ such that
$|x|=\ff(|x_0|,\,|x_1|)$, where $x_0\in X_0,\quad x_1\in X_1$ with the norm
$$
\|x\|_{\ff(X_0,X_1)}=\inf_{x_0,x_1}\max (\|x_0\|_{X_0},\,\|x_1\|_{X_1})
$$
and the infimum is taken over all representation of $|x|$ in the form
$\ff(|x_0|,|x_1|)$.

\smallskip
\noindent{\bf Lemma 1.} \textsl{Let $1\le p_0,p_1\le\ii$, then}
$$
\ff(\Lppp)=\ff(\Lppp)_{p_0,p_1}.
$$
(Note that if $U_0=1$ and $U_1=1$, then $\ff(L_{p_0},\,L_{p_1})$
is an Orlicz space.)

\smallskip
\noindent{\it Proof.} Let us show the embedding
\begin{equation}\label{star}
\ff(\Lppp)_{p_0,p_1}\subset \ff(\Lppp).
\end{equation}

Let
$$
x=\sum_{n\in \Bbb Z}w_n\rr(t_n)
$$
where $\{\|w_n\|_{L_{p_0}}\}\in l_{p_0}$, and
$\{t_n\|w_n\|_{L_{p_1}}\}\in l_{p_1}$. We introduce a linear
operator
$$
T:\{l_{p'_0},\,l_{p'_1}(t_n^{-1})\}\to\Lpp,
$$
where $p'$ denotes the index dual to the index $p$, i.e.,
$p'=(1-1/p)^{-1}$. We put
$$
T(\xi)=\sum_{n\in \Bbb Z}w_n\xi_n,\qquad({\rm convergence\ in }\
L_{p_0}(U_0)+L_{p_1}(U_1))
$$
for $\xi\in l_{p'_0}+l_{p'_1}(t_n^{-1})$.

If $\xi\in l_{p'_0}$, then the series $(w_n\xi_n)$ is absolutely
convergent in $L_{p_0}(U_0)$. If $\xi\in l_{p'_1}(\tttt)$,
then the series $(w_n\xi_n)=(w_n\ttt\xi_n\tttt)$ is absolutely
convergent in $L_{p_1}(U_1)$. Moreover everywhere on the measure space
$\frak M$ we have
$$
|\sum_{n\in\Bbb Z}w_n({m})\xi_n|\le\sum_{n\in\Bbb Z}
|w_n({m})||\xi_n|
$$
$$
\le
(\sum_{n\in\Bbb Z}|w_n({m})|^{p_0})^{1/p_{0}}
(\sum_{n\in\Bbb Z}|\xi_n|^{p'_0})^{1/p'_0}.
$$

Since
$$
\int\limits_{\frak M}\sum_{n\in\Bbb Z}|U_0({m})w_n({m})|^{p_0}d{m}
=\sum_{n\in\Bbb Z}\|w_n\|^{p_0}_{L_{p_0}(U_0)}<\ii
$$
we deduce
$$
(\sum_{n\in\Bbb Z}|w_n({m})|^{p_0})^{1/p_0}=W_0({m})\in
L_{p_0}(U_0).
$$

Analogously
$$
(\sum_{n\in\Bbb Z}(\ttt |w_n({m})|)^{p_1})^{1/p_1}=W_1({m})\in
L_{p_1}(U_1).
$$

Hence
$$
|T(\xi)({m})|\le W_0({m})\|\xi\|_{l_{p'_0}} \quad{\rm for}\quad
\xi\in l_{p'_0},
$$
$$
|T(\xi)({m})|\le W_1({m})\|\xi\|_{l_{p'_1}(\tttt)}
 \quad{\rm for}\quad\xi\in l_{p'_1}(\tttt).
$$
This means that the operator $T$ is factorized through a couple of
weighted $L_\ii$-spaces. So we have
$$
T:l_{p'_0}\to L_\ii(W_0^{-1}({m}))\subset L_{p_0}(U_0),
$$
$$
T:l_{p'_1}(\tttt)\to L_\ii(W_1^{-1}({m}))\subset L_{p_1}(U_1).
$$
(Without loss of generality we suppose here that almost everywhere
on the measure space $\frak M$ both weights $W_0$ and $W_1$ are
not equal to 0.)

It is easily seen that the image $x$ of the sequence $\{\rr(\ttt)\}$
with respect to the operator $T$ is equal to the sum
$$
x=\sum_{n\in\Bbb Z}\rr(\ttt)w_n.
$$

In view of (\ref{J})
$$
K(s,t,\{\rr(\ttt)\},\{l_{p'_0},l_{p'_1}(\tttt)\})\asymp\ff(s,t).
$$
Hence
$$
K(s,t,T(\{\rr(\ttt)\}),\{L_{\ii}(W_0^{-1}),L_{\ii}(W_1^{-1})\})
$$
$$=
K(s,t,x,\{L_{\ii}(W_0^{-1}),L_{\ii}(W_1^{-1})\})\le C\ff(s,t),
$$
which imply
$$
|x({m})|\le
K(W_0({m}),W_1({m}),x,\{L_{\ii}(W_0^{-1}),L_{\ii}(W_1^{-1})\})
\le C\ff(W_0({m}),W_1({m})).
$$
So we obtain $x\in \ff(\Lppp)$ and embedding (\ref{star}) is proved.

Now consider the embedding
\begin{equation}\label{2star}
\ff(\Lppp)\subset\ff(\Lppp)_{p_0,p_1}.
\end{equation}

Let $a\in \ff(\Lppp)$. Denote by $\psi(t)$ the K-functional of $a$, i.e.,
$$
\psi(t)=K(t,a,\Lpp).
$$
and denote by $\{\uu\}$ a balanced sequence for the function
$\pp$. (Note that $\pp\in\Phi_0$ because of $a\in
\ff(\Lppp)\subset (\Lppp)_{\rr,\ii}$.) Then
$$
\pp(t)\asymp K(t,\{\pp(\uu)\},\{\lpu\})\asymp
K(t,\{\pp(\uu)\},\{l_\ii,l_\ii(\uuu)\}).
$$
(Note also that index $m$ runs over some interval $\Bbb M$ in $\Bbb Z$.)
By the Sparr theorem for the couples $\Lpp$ and $\{\lpu\}$ there exist
linear operators
$$
A:\Lpp\to\{\lpu\}
$$
and
$$
B:\{\lpu\}\to\Lpp
$$
such that $A(a)=\{\pp(\uu)\}$ and $B(\{\pp(\uu)\})=a$. Since
the Calderon--Lozanovskii construction is an interpolation functor
for couples of $L_p$ spaces we deduce
$$
\{\pp(\uu)\}\in \ff(\lpu).
$$

Recall that any element $a\in\ff(\Lppp)$ belongs to the orbit of the sequence
$\{\rr(\ttt)\}$ with respect to linear operators mapping the couple
$\{l_\ii,l_\ii(\tttt)\}$ into the couple $\Lpp$
(see, for instance \cite{8}). This means that there exists a linear operator
$$
T:\{l_\ii,l_\ii(\tttt)\}\to \Lpp
$$ such that $T(\{\rr(\ttt)\})=a$.

Let us consider the composition
$$
AT:\{l_\ii,l_\ii(\tttt)\}\to \{\lpu\},
$$
which maps $\{\rr(\ttt)\}$ into $\{\pp(\uu)\}$.

Since $l_p\subset l_\ii$ we consider the operator $AT$ as an
operator mapping the couple $\{l_\ii,l_\ii(\tttt)\}$ into the
couple $\{l_\ii,l_\ii(\uuu)\}$. The embedding $l_p\subset l_\ii$
are $(1,p)$--summing by the Karl--Bennett theorem \cite{9}, hence
$$
\{\|AT(e_n)\|_{l_\ii}\}\in l_{p_0}\qquad
\{\ttt\|AT(e_n)\|_{l_\ii(\uuu)}\}\in l_{p_1},
$$
where $\{e_n\}$ denotes the standard basis in $l_p$ spaces.

Hence
\begin{equation}
\label{l2}
{\pp(\uu)}=\sum_{n\in\Bbb Z}\rr(\ttt)AT(e_n)=
\sum_{n\in\Bbb Z}\rr(\ttt)w_n,
\end{equation}
which means
$$
\{\pp(\uu)\}\in \ff(l_\ii,l_\ii(\uuu))_{p_0,p_1}.
$$

In view of (\ref {l2})
$$
K(t,\{\pp(\uu)\}, \{l_\ii,l_\ii(\uuu)\})\le
\sum_{n\in\Bbb Z}\rr(\ttt)
K(t,w_n, \{l_\ii,l_\ii(\uuu)\}),
$$
therefore by (\ref{J})

\begin{equation}\label{l3}
K(t,\{\pp(\uu)\}, \{\lpu\})\le
C\sum_{n\in\Bbb Z}\rr(\ttt)
K(t,w_n, \{l_\ii,l_\ii(\uuu)\}).
\end{equation}

By the K-divisibility (see \cite {BK}) for the couple $\{\lpu\}$ (\ref {l3})
implies that there exists a sequence $w'_n\in l_{p_0}\cap l_{p_1}(\uuu)$
such that
$$
\{\pp(\uu)\}=\sum_{n\in\Bbb Z}\rr(\ttt)w'_n,
$$
and
\begin{equation}\label{l4}
K(t,w'_n, \{\lpu\})\le 6CK(t,w_n, \{l_\ii,l_\ii(\uuu)\}).
\end{equation}
Since
$$
\lim_{t\to\ii}K(t,w'_n,\{\lpu\})=\|w'_n\|_{l_{p_0}}
$$
and
$$
\lim_{t\to 0}\frac1tK(t,w'_n,\{\lpu\})=\|w'_n\|_{l_{p_1}(\uuu)}
$$
inequalities (\ref {l4}) imply
$$
\|w'_n\|_{l_{p_0}}\le 6C\|w'_n\|_{l_\ii},\qquad
\|w'_n\|_{l_{p_1}(\uuu)}\le 6C\|w'_n\|_{l_\ii(\uuu)}.
$$

Hence
$$
\{\pp(\uu)\}\in \ff(\lpu)_{p_0,p_1}.
$$

The construction $\ff(X_0,X_1)_{p_0,p_1 }$ is evidently an interpolation
functor, hence the image of $\{\pp(\uu)\}$ with respect to
linear operator $B:\{\lpu\}\to\Lpp$ belongs to $\ff(\Lppp)_{p_0,p_1}$.
Namely $a=B(\{\pp(\uu)\})\in\ff(\Lppp)_{p_0,p_1}$. So the embedding
(\ref{2star}) and Lemma 1 are proved.

{\bf Remark 1.} The embedding (\ref{2star})  in the case $p_0,\,p_1\ge 2$
may be proved without the Karl--Bennett theorem and K-divisibility.
We may simply mention that for a unconditionally convergent series
$(w_n)$ in $L_p$ we have $\{\|w_n\|_{L_p}\}\in l_p$
if $p\ge 2$.

Let us now describe the space $\ff(\Lppp)$ in term of the K-functional.
Take any $a\in\ff(\Lppp)$, consider the K-functional
$\psi(t)= K(1,t,a,\{\Lppp\})$ of $a$.
Since $\ff\in\Phi_0$, then $\psi\in\Phi_0$.

It is easily seen that $a\in \ff(\Lppp)$ if and only if
$\{\psi(\uu)\}\in \ff(l_{p_0},\,l_{p_1}(\uuu))$, where $\{\uu\}$
is a balanced sequence for the function.

Recall that interpolation function $\ff$ is called non-degenerate if the ranges
of the functions $\ff(t,1)$ and $\ff(1,t)$ where $t>0$ coincide with
$(0,\ii)$.

\smallskip
\noindent{\bf Lemma 2.} \textsl{If $\ff$ is non-degenerate, the
space $\ff(X_0,\,X_1)_{p_0,p_1}$ consists of $x\in X_0+X_1$ for
which
$$
\{K(\uu,x,\XX\}\in \ff(l_{p_0},\,l_{p_1}(\uuu)),
$$
where $\{\uu\}$ is a balanced sequence for the function
$K(t,x,\{X_0,X_1\})$.}

\smallskip
{\bf Remark 2.} Note that the spaces $\ff(X_0,\,X_1)_{p,p}$
coincide with the space $(X_0,\,X_1)_{\rr,p}$ introduced by
S.Janson (see \cite{10}). Lemma 2 gives us a new description of
these spaces as well.

\smallskip
\noindent {\it Proof.} Let us show that Lemma 2 is true for any interpolation
function $\ff\in\Phi_0$ and any mutually closed Banach couple.
Recall that a  couple $\{X_0,X_1\}$ is called mutually closed if
$X_0=\widetilde X_0$ and $X_1=\widetilde X_1$, where
$\widetilde X$ denotes the completion of $X$ with respect to
$X_0+X_1$.

If $x\in \ff(X_0,X_1)_{p_0,p_1}$, then
\begin{equation}
\label{3star}
x=\sum_{n}\rr(\ttt)w_n
\end{equation}
where
$$
\sum_{n}\|w_n\|^{p_0}_{X_0}<\ii
,\qquad\sum_{n}(\ttt\|w_n\|_{X_1})^{p_1}<\ii.
$$

The expansion (\ref{3star}) implies
$$
K(t,x,\XX)\le \sum_{n}\rr(\ttt)K(t,w_n,\XX).
$$
Let us denote again by $\pp(t)=K(t,x,\XX)$. Hence
$$
K(t,\{\pp(\uu)\}, \{\lpu\})\le
\sum_{n}\rr(\ttt)K(t,w_n,\XX).
$$
By the K-divisibility applied to the couple $\{\lpu\}$ we obtain
a sequence $w'_n\in l_{p_0}\cap l_{p_1}(\uuu)$ such that
$$
\{\pp(\uu)\}=\sum_{n}\rr(\ttt)w'_n
$$
and
\begin{equation}\label{4star}
K(t,w'_n,\{\lpu\})\le 6C K(t,w_n,\XX)
\end{equation}
for any $n$.

By (\ref{4star}) we get
$$
\|w'_n\|_{l_{p_0}}\le6C\|w_n\|_{X_0},
\qquad
\|w'_n\|_{l_{p_1}(\uuu)}\le6C\|w_n\|_{X_0}
$$
as it was done in Lemma 1. So we conclude $\{\pp(\uu)\}\in\ff(\lpu)_{p_0,p_1}$
and
$$
\{\pp(\uu)\}\in \ff(\lpu)
$$
by Lemma 1.

Now we suppose that
$$
\{\pp(\uu)\}\in \ff(\lpu)=\ff(\lpu)_{p_0,p_1},
$$
then
$$
\{\pp(\uu)\}=\sum_{n}\rr(\ttt)w'_n,
$$
where
$$
\{\|w'_n\|_{l_{p_0}}\}\in l_{p_0},
\qquad\{\|w'_n\|_{l_{p_1}(\uuu)}\}\in l_{p_1}.
$$

Hence
$$
K(t,a,\XX)=\pp(t)\le CK(t,\{\pp(\uu)\},\{\lpu\})
$$
$$
\le C\sum_{n}\rr(\ttt)K(t,w_n,\{\lpu\}).
$$

By the K-divisibility of the couple $\XX$ we find a sequence $w_n\in X_0+X_1$
such that
$$
a=\sum_{n}\rr(\ttt)w_n
$$
and
$$K(t,w_n,\XX)\le 6CK(t,w'_n,\{\lpu\}).
$$

Hence
$$
\lim_{t\to \ii}K(t,w_n,\XX)\le6C\|w'_n\|_{l_{p_0}},
$$
$$
\lim_{t\to 0}\frac1tK(t,w_n,\XX)\le6C\|w'_n\|_{l_{p_1}(\uuu)}.
$$
This means
$$
\|w_n\|_{\widetilde X_0}\le 6C\|w'_n\|_{l_{p_0}},\qquad
\|w_n\|_{\widetilde X_1}\le 6C\|w'_n\|_{l_{p_1}(\uuu)}.
$$

Thus  $a\in \ff(\widetilde X_0,\widetilde X_1)_{p_0,p_1}$.
So if $X_0=\widetilde X_0$ and $X_1=\widetilde X_1$, Lemma 2 is proved
without any assumption on the function $\ff\in\Phi_0$. It may be
shown that
$$
\ff(X_0,X_1)_{p_0,p_1}=\ff(\widetilde X_0,\widetilde X_1)_{p_0,p_1}
$$
if $\ff$ is non-degenerate (see \cite{K}). Lemma 2 is proved.

Let now $\ff\in\Phi_0$ and $\ff$ be degenerate. If both $\ff(1,t)$
and $\ff(t,1)$ are bounded then $\ff(t,s)\asymp\min(s,t)$ and naturally
we have
$$
\ff(X_0,X_1)_{p_0,p_1}=X_0\cap X_1.
$$

Suppose now that $\ff(1,t)$ is bounded and $\ff(t,1)$ is not
bounded. In this case a balanced sequence for $\ff$ turns out to
be one-sided. Let us think that we have $\{\ttt\}^0_{-\ii}$.

Let us show that
$$
\ff(X_0,X_1)_{p_0,p_1}=(X_0,X_1)^K_{\ff,p_0,p_1}\cap X_0,
$$
where by $(X_0,X_1)^K_{\ff,p_0,p_1}$ we denote the space of $x\in X_0+X_1$
such that
$$
\{K(\uu,x,\XX)\}\in \ff(\lpu).
$$

We have already shown that
$$
(X_0,X_1)^K_{\ff,p_0,p_1}=\ff({\widetilde X_0},{\widetilde X_1})_{p_0,p_1}.
$$

Since $\{\ttt\}$ is one-sided we have
$$
\ff(X_0,X_1)_{p_0,p_1}\subset\ff(X_0,X_1)_{\ii,\ii}\subset X_0.
$$

Hence
$$
\ff(X_0,X_1)_{p_0,p_1}\subset(X_0,X_1)^K_{\ff,p_0,p_1}\cap X_0.
$$

If now $x\in (X_0,X_1)^K_{\ff,p_0,p_1}\cap X_0$, then there exists
an element $a\in (l_1,l_1(2^{-k}))_{\ff,p_0,p_1}$ such that
$x=Ua$, where $U:\{l_1,l_1(2^{-k})\}\to\XX$ (see \cite {8}, where
we have shown that interpolation from the couple  $\{l_1,l_1(2^{-k})\}$ to
any Banach couple is essentially described by K-method.) The couple
$\{l_1,l_1(2^{-k})\}$  is mutually closed. Hence
$a \in \ff(l_1,l_1(2^{-k}))_{p_0,p_1}$, therefore
$$
a=\sum_{n=-\ii}^0 \rr(\ttt)a_n,
$$
where $\{\|a_n \|_{l_1}\}\in l_{p_0}$,
$\{\ttt\|a_n \|_{l_1(2^{-k})}\}\in l_{p_1}$. Thus
$$
x=Ua=\sum_{n=-\ii}^0 \rr(\ttt)U(a_n),
$$
where $\{\|U(a_n)\|_{X_0}\}\in l_{p_0}$,
$\{\ttt\|U(a_n)\|_{X_1}\}\in l_{p_1}$. So we have
$x\in \ff(X_0,X_1)_{p_0,p_1}$.

Thus we obtain a description of the functor
$\ff(X_0,X_1)_{p_0,p_1}$ for all functions $\ff\in \Phi_0$.

If at least one of $p_0,\,p_1$ is equal to infinity, then
functions $K(s,t,a,\Lpp)$ where $a\in L_{p_0}(U_0)+L_{p_1}(U_1)$
not always belong to $\Phi_0$. So if we are going to describe
orbits in terms of $\ff(X_0,X_1)_{p_0,p_1}$, we need in the
definition of the functor $\ff(X_0,X_1)_{p_0,p_1}$ not only for
$\ff\in \Phi_0$.

Recall that any interpolation function $\ff(s,t)$ can be presented
in the form
$$
\ff(s,t)=\alpha s+\beta t +\ff_0(s,t),
$$
where $\ff_0 \in \Phi_0$. If $\alpha$ or $\beta$ is greater than $0$, then
$\ff_0 \notin \Phi_0$.

\smallskip
\noindent{\bf Definition 2.} \textit{If $\ff$ is arbitrary
interpolation function, then denote by $\ff(X_0,X_1)_{p_0,p_1}$
the space
$$
\alpha X_0+\beta X_1 +\ff_0(X_0,X_1)_{p_0,p_1},
$$
where $\alpha X$ denotes the space $X$ with the norm
$\|x\|_{\alpha X}=\alpha \|x\|_X$ if $\alpha >0$, and $\alpha X=0$
if $\alpha =0$.}

\smallskip
Note that if both $\alpha$ and $\beta$ are greater than zero, then
$$
\ff(X_0,X_1)_{p_0,p_1}=X_0+X_1
$$
for any $p_0$ and $p_1$.

$2^o$. {\sc The main theorem}.

\noindent{\bf Theorem.} \textsl{Let $\{\Lppp\}$ and $\{\Lqqq\}$ be
two Banach couples, where $1\le p_0,p_1,q_0,q_1\le\ii$, and $a\in
L_{p_0}(U_0)+L_{p_1}(U_1)$ such that
$\ff(s,t)=K(s,t,a,\,\{\Lppp\})\in \Phi_0$, then
$$
Orb(a,\{\Lppp\}\to\Lqq)=\ff(\Lqqq)_{r_0,r_1},
$$
where $\0$ and $\1$. (As usual $x_+$ denotes the positive part of
$x$.)}

The rest cases $\ff(s,t)\notin \Phi_0$ can be easily reduced to
$\ff(s,t)\in \Phi_0$ as it was done in \cite {8} where the analogous
situation takes place for $p_0 \le q_0$ and $p_1 \le q_1$.

The proof is a combination of the following propositions.

\smallskip
\noindent{\bf Proposition 1.} \textsl{For any $1 \le
p_0,p_1,q_0,q_1 \le\ii$ and any weights $U_0,\,U_1,\,V_0,\,V_1$,
and for any $a\in L_{p_0}(U_0)+L_{p_1}(U_1)$}
$$
Orb(a,\{\Lppp\}\to\Lqq)\subset\ff(\Lqqq)_{r_0,r_1}.
$$

\smallskip
\noindent \textit{Proof.} Let $b=Ta$, where $T:\Lpp\to\Lqq$.
Recall that $\rr(t)=\ff(1,t)$. Denote $a_{\rr}=\{\rr(\ttt)\}$,
$\psi(u)=K(u,b,\Lqq)$ and $b_\psi=\{\pp(u_n)\}$, where $\uu$ is a
balanced sequence for $\psi(u)$. The Sparr theorem implies that
there exists a linear operator $S:\lpp\to\lqq$ such that
$Sa_{\rr}=b_{\psi}$.

We consider the embedding $\lqq\subset\{l_\ii,\,l_\ii(\uuu)\}$. It
is known that the embedding $l_{q_i}\subset\l_\ii$ are $(1,q_i)$--
summing operators (by the Karl--Bennett theorem, see \cite{9}).
Hence if $q _0<p_0$, then the image of the standard basic sequence
in $l_{p_0}$ with respect to $S:l_{p_0}\to l_{\ii}$ is
$l_{r_0}$-sequence, that is $ \{\|S(e_n)\|_{l_\ii}\}\in l_{r_0}, $
where $r_0^{-1}=q_0^{-1}-p_0^{-1}$. Analogously $
\{t_n\|S(e_n)\|_{l_\ii(\uuu)}\}\in l_{r_1}, $ where
$r_1^{-1}=q_1^{-1}-p_1^{-1}$. Hence in any case we have
$$
\{\|S(e_n)\|_{l_\ii}\}\in l_{r_0},\qquad {\rm and}
\qquad \{t_n\|S(e_n)\|_{l_\ii(\uuu)}\}\in l_{r_1},
$$
where $\0$ and $\1 $.

Therefore by definition $b_\pp\in \ff(l_\ii,\,l_\ii({\uuu}))_{r_0,r_1}$.
By Lemma 2 this means
$$
\{K(v_m,b_\pp,\,\{l_\ii,l_\ii(u_m^{-1})\})\}\in
\ff(l_{r_0},\,l_{r_1}(v_m^{-1})),
$$
where $\{v_m\}$ is a balanced sequence of the function
$K(v,b_\pp,\,\{l_\ii,l_\ii(\uuu)\})\asymp
K(v,b_\pp,\,\lqq)=\pp(v)$.
 Hence we can take $v_m=\uu$, and by (\ref {J})
 $$
K(\uu,b,\Lqq)\asymp K(\uu,b_\pp,\lqq)
\asymp K(\uu,b_\pp,\,\{l_\ii,l_\ii(\uuu)\})
$$

So $\{K(\uu,b,\Lqq)\}\in\ff(\lrrr)$. By Lemma 2 Proposition is proved.

The following propositions are devoted to the inverse inclusion
$$
\ff(\Lqqq)_{r_0,r_1}\subset Orb(a,\Lpp\to\Lqq).
$$

For any $b\in \pp(\Lqqq)_{r_0,r_1}$ we must find an operator
$T\in \Lpp\to\Lqq$ such that $b=Ta$. Again with the help of the Sparr
theorem we substitute $a$ by $a_\rr$ and $b$ by $b_\pp$ as well as
initial couples by $\lpp $ and $\lqq$ respectively.

\smallskip
\noindent{\bf Proposition 2.} \textsl{Let
$\{\pp(\uu)\}\in\ff(\lrrr)$, then there exist sequences
$\{\beta_m^0\}\in l_{r_0}$ and $\{\beta_m^1\}\in l_{r_1}$ such
that}
$$
K(s,t,\{\pp(\uu)\},\,\{l_1(1/\beta_m^0),\,l_1(1/\beta_m^1\uu)\})\le
C\,\ff(s,t).
$$

\smallskip
\noindent{\it Proof.}
Since $\{\pp(\uu)\}\in \ff(\lrrr)$ we have $\pp(\uu)=
\ff(\alpha_m^0,\alpha_m^1\uu)$, where
$\{\alpha_m^0\}\in l_{r_0}$, $\{\alpha_m^1 \}\in l_{r_1}$.

Consider now the convolution of $\alpha_m^0$ and  $\varepsilon_j^0$,
where  $\varepsilon_j^0=(1-\varepsilon)^{-j}$ for $j \le 0$ and
$\varepsilon_j^0=0$ for $j>0$. Denote $\beta^0=\varepsilon^0*\alpha^0$.
Analogously $\beta^1=\varepsilon^1*\alpha^1$, where
$\varepsilon_j^1=(1-\varepsilon)^{j}$ for $0 \le j$ and
$\varepsilon_j^0=0$ for $j<0$. Choose  $\varepsilon>0$ sufficiently small.
Then
\begin{equation}\label{ee}
\beta^0_{m+1}(1-\varepsilon)\le\beta^0_{m},\qquad
\beta^1_{m}(1-\varepsilon)\le\beta^1_{m+1},
\end{equation}
$\alpha^0_m \le \beta^0_m$, $\alpha^1_m \le \beta^1_m$
and $\beta^0 \in l_{r_0}$, $\beta^1 \in l_{r_1}$.

Indeed $\varepsilon^0,\ \varepsilon^1\in l_1$, hence
$\beta^0\in l_{r_0}$, $\beta^1 \in l_{r_1}$ in view of convolution
properties. Since $\varepsilon^0_0=\varepsilon^1_0=1$ we have
$\alpha^0_m\le \beta^0_m$, $\alpha^1_m\le \beta^1_m$.
By definition
$$
\beta^0_m=\sum_{k=-\ii}^0(1-\varepsilon)^k\alpha^0_{m-k}.
$$
(We define $\alpha^0$ as 0 outside the interval $\Bbb M$.)

Since
$$
\beta^0_{m+1}=\sum_{k=-\ii}^{-1}(1-\varepsilon)^{k+1}\alpha^0_{m-k}=
(1-\varepsilon)\sum_{k=-\ii}^{-1}(1-\varepsilon)^{k}\alpha^0_{m-k}
$$
we obtain
$\beta^0_{m+1}(1-\varepsilon)\le\beta^0_{m}$.

Analogously
$$
\beta^1_m=\sum_{k=0}^{\ii}(1-\varepsilon)^k\alpha^1_{m-k},
$$
and
$$
\beta^1_{m+1}=\sum_{k=-1}^{\ii}(1-\varepsilon)^{k+1}\alpha^1_{m-k}=
(1-\varepsilon)
\sum_{k=-1}^{\ii}(1-\varepsilon)^{k}\alpha^1_{m-k}.
$$
Hence
$\beta^1_{m}(1-\varepsilon)\le\beta^1_{m+1}$.

For the K-functional of $\{\pp(\uu)\}$ we have now
$$
K(s,t,\{\pp(\uu)\},\,\{l_1(1/\beta_m^0),\,l_1(1/\beta_m^1\uu)\})=
\sum_{m\in \Bbb M}\pp(\uu)
\min\left( \frac{s}{\beta^0_m},\,\frac{t}{\beta^1_m\uu}
\right).
$$

For $s=\beta^0_k$, $t=\beta^1_k u_k$ we have
$$
K(\beta^0_k,\,\beta^1_k u_k)\le
\sum_{m\le k}\pp(\uu)\frac{\beta^0_k}{\beta^0_m}+
\sum_{m\ge k}\pp(\uu)\frac{\beta^1_k u_k}{\beta^1_m\uu}.
$$
By (\ref{J_n}) and (\ref{ee}) we conclude that right hand side of the latter
inequality is estimated by a sequence equivalent to $\pp(u_k)$.

Indeed in view of (\ref{J_n})
$$
\frac{\pp(\uu)}{\pp(u_{m+1})}\le\frac 1q,\qquad
\frac{\uu}{\pp(\uu)}\frac{\pp(u_{m+1})}{u_{m+1}}\le\frac 1q,
$$
and by (\ref{ee})
$$
\frac{\beta^0_{m+1}}{\beta^0_{m}}\le\frac 1{1-\varepsilon},
\qquad
\frac{\beta^1_{m}}{\beta^1_{m+1}}\le\frac 1{1-\varepsilon}.
$$
Therefore
$$
\frac{\pp(\uu)}{\beta^0_{m}}
\frac{\beta^0_{m+1}}{\pp(u_{m+1})}\le\frac 1{q'},
\qquad
\frac{{\beta^1_{m}}\uu}{\pp(\uu)}\frac{\pp(u_{m+1})}{{\beta^1_{m+1}}u_{m+1}}
\le\frac 1{q'}.
$$
where $q'>1$ if $\varepsilon$ is sufficiently small.

Hence
$$
\frac{\pp(\uu)}{\beta^0_{m}}
\frac{\beta^0_{k}}{\pp(u_{k})}\le\left(\frac 1{q'}\right)^{k-m},
$$
for $m\le k$, and
$$
\frac{{\beta^1_{m}}\uu}{\pp(\uu)}\frac{\pp(u_{k})}{{\beta^1_{k}}u_{k}}
\le\left(\frac 1{q'}\right)^{m-k}
$$
for $m\ge k$. These inequalities imply
$$
\sum_{m\le k}\pp(\uu)\frac{\beta^0_k}{\beta^0_m}\le
\frac{q'}{q'-1}\pp(u_k),
\qquad
\sum_{m\ge k}\pp(\uu)\frac{\beta^1_k u_k}{\beta^1_m\uu}\le
\frac{q'}{q'-1}\pp(u_k).
$$

Thus
$$
K(\beta^0_k,\,\beta^1_k u_k)\le
\frac{q'}{q'-1}\,\pp(u_k)=\frac{q'}{q'-1}\,\ff(\alpha_k^0,\alpha_k^1u_k)\le
\frac{q'}{q'-1}\,\ff(\beta_k^0,\beta_k^1u_k),
$$
therefore
$$
K(s,t,\{\pp(\uu)\},\,\{l_1(1/\beta_m^0),\,l_1(1/\beta_m^1\uu)\})\le
C\,\ff(s,t).
$$
Proposition is proved.

\smallskip
\noindent{\bf Proposition 3.} \textsl{Let
$b_\pp=\{\pp(\uu)\}\in\ff(\lrrr)$, then there exists a linear
operator $S:\lpp\to\lqq$ such that $S(a_\rr)=b_\pp$.}

\smallskip
\noindent{\it Proof.} Without loss of generality we assume that
$p_0\ge q_0$, $p_1\ge q_1$. By Proposition 2 we can find
$\beta^0\in l_{r_0}$ and $\beta^1\in l_{r_1}$. Consider embedding
\begin{equation}\label{inc}
\{l_1(1/\beta_m^0),\,l_1(1/\beta_m^1\uu)\}\subset
\{l_{p_0}(1/\beta_m^0),\,l_{p_1}(1/\beta_m^1 \uu)\subset\lqq\}
\end{equation}
and the element $b_\pp$.

By the left-hand side embedding and Proposition 2 we have
$$
K(s,t,\{\pp(\uu)\},\,\{l_{p_0}(1/\beta_m^0),\,l_{p_1}(1/\beta_m^1\uu)\})\le
C\,\ff(s,t).
$$
Since $\ff(s,t)\asymp K(s,t,a_\rr,\lpp)$ by the Sparr theorem there exists
an operator $S:\lpp\to \{l_{p_0}(1/\beta_m^0),\,l_{p_1}(1/\beta_m^1 \uu)\}$
mapping $a_\rr$ into $b_\pp$.

The composition of $S$ and the right hand side embedding in (\ref{inc})
is the desired mapping. Thus Proposition and Theorem are proved.

 \end{document}